\documentclass{amsart}

\usepackage{amsrefs}

\usepackage{amsmath}
\usepackage{amssymb}
\usepackage{enumerate}

\theoremstyle{plain} %% This is the default

\newtheorem{thm}{Theorem}

\theoremstyle{definition}

\theoremstyle{remark}

\numberwithin{equation}{section}
\newcommand{\thmref}[1]{Theo\-rem \ref{#1}}

\title{Optimal algorithms of Gram-Schmidt type}
\author{James B. Wilson}
\address{
	Department of Mathematics\\
	The Ohio State University\\
	Columbus, OH 43210\\
}
\email{wilson@math.ohio-state.edu}
\date{\today}
\begin{document}

\maketitle

\begin{abstract}
Three algorithms of Gram-Schmidt type are given that produce an orthogonal decomposition of finite $d$-dimensional symmetric, alternating, or Hermitian forms over division rings.  The first uses $d^3/3+O(d^2)$ ring operations with very simple implementation.  Next, that algorithm is adapted in two new directions.  One is an optimal sequential algorithm whose complexity matches the complexity of matrix multiplication.  The other is a parallel NC algorithm with similar complexity.
\end{abstract}

%--------------------
\section{Introduction}

The classic Gram-Schmidt `orthogonalization process' returns an orthonormal basis of an inner product space.  Here we generalize that process in the appropriate fashion to Hermitian forms over division rings $\Delta$.  For us a Hermitian $\Delta$-form is a function $b:V\times V\to \Delta$ on a finite-dimensional $\Delta$-vector space $V$ where $b$ is linear in the first variable and for some anti-isomorphism $\sigma$ of $\Delta$, for all $u,v\in V$, $b(u,v)=b(v,u)^{\sigma}$.  This captures the usual symmetric and skew-symmetric forms as well as the traditional Hermitian forms; cf \cite{Taylor}.   We identify $V$ with a space of row vectors and so describe $b$ by a matrix $B$ such that $b(u,v)=uBv^{\sigma t}$.  The assumptions on $b$ force $B=0$, or $B=s B^{\sigma t}$ with $s=\pm 1$ and $\sigma^2=1$.  To change the basis we use an invertible matrix $A$ and observe $b(uA,vA)=uABA^{\sigma t} v^{\sigma t}$.  Hence, a fully refined orthogonal decomposition for $b$ is captured by a matrix $A$ under which $ABA^{\sigma t}$ is nearly diagonal, nearly in that sometimes $J:=\begin{bmatrix} 0 & 1 \\ -1 & 0 \end{bmatrix}$ is required.

\begin{thm}\label{thm:main}
Let $\Delta$ be a division ring, let $s=\pm 1$, and let $\sigma$ be 
a unital anti-isomorphism of $\Delta$ with $\sigma^2=1$.
There are deterministic algorithms that, given a $(d\times d)$-matrix $B=sB^{\sigma t}$, return an invertible $(d\times d)$-matrix $A$ such that
\begin{equation}\label{eq:standard}
	ABA^{\sigma t}  = B_1\oplus \cdots \oplus B_m
\end{equation}
and each $B_i$ is either $1\times 1$ or $J$.
%\begin{equation}\label{eq:standard}
%\begin{split}
%	ABA^{\sigma t} 
%		& = \begin{bmatrix} 0 & 1 \\ s & 0\end{bmatrix}\oplus \cdots\oplus
%			\begin{bmatrix} 0 & 1 \\ s & 0\end{bmatrix}\oplus [\beta_{1}]\oplus \cdots 
%				\oplus [\beta_{m}].
%\end{split}
%\end{equation}
\begin{enumerate}[(i)]
\item
The first algorithm uses $e$ inversions, $\binom{d}{2}$ equality tests, $e^3/3+O(d^2)$ additions and $e^3/3+O(d^2)$ multiplications in $\Delta$, where $e=d-r$ and $r$ is the rank of the null space of $B$.

\item
The second algorithm returns a straight line program to $A$ using $O(d^{\omega})$ operations in $\Delta$, where $\omega$ is the exponent of matrix multiplication.

\item
The third algorithm is parallel $NC^3$ in an arithmetic model, that is, it uses $O(\log^3 d)$ operations in $\Delta$ on $d^{O(1)}$ processors.
\end{enumerate}
\end{thm}

\thmref{thm:main} in part proves that Hermitian forms over division rings have a decomposition 
of the type described in \eqref{eq:standard}.  For fields this is well-known, e.g. \cite[Theorems 3.7]{Artin:geom}, but most proofs begin by classifying forms into subclasses, e.g. symmetric if $s=1$ and $\sigma=1$, skew-symmetric if $s=-1$ and $\sigma=1$, or Hermitian if $\sigma\neq 1$.    Nice bases are constructed by individual arguments for each case.  Here we find a single argument allows for uniform optimal asymptotic and parallel algorithms without dependence on $\Delta$.  

The idea behind \thmref{thm:main}(i) is shared by many generalizations of Gram-Schmidt.  For symmetric forms it goes back at least to Smiley's \emph{Algebra of Matrices} \cite[Section 12.2]{Smiley} and is adequately described as symmetric Gaussian elimination.  Dax and Kaniel \cite{Dax-Kaniel} give a detailed analysis of such an algorithm for symmetric forms. Holt and Roney-Dougal \cite{HRD} use the method in a case-by-case algorithm for Hermitian forms over finite fields.  I was also gratefully alerted to a predecessor to \thmref{thm:main}(ii) that applies to symmetric forms over fields; see \cite[Theorem 16.25]{Burgisser:comp}.

The algorithms for \thmref{thm:main} parts (ii) and (iii) settle the complexity of finding an nice basis for a generic Hermitian form, but these may not be best suited for certain applications.  First, they depend on data structures for fast matrix multiplication which may provide an undesirable overhead in small dimensions.  The exact cross-over dimension is an issue of ongoing research; see \cite[p. 313]{vzG}.  Furthermore, our algorithms assume exact field operations, such as in algebraic number fields, rational Quaternion division rings, or  finite fields. We make no claims about their numerical stability in fields with floating point approximations.  In such cases consider \cite{Lingen}.

Each of our algorithms allows the user to choose a computational encoding for $\Delta$, such as by polynomials or matrices over a field.  If no alternative suggests itself, an adequate method is to encode $\Delta$ by structure constants over its center; see \cite[p. 223]{Ronyai}.
Also, $b$ can be encoded as a ``black-box''; however, we will eventually evaluate $b$ on all unordered pairs from a fixed basis for $V$ and so it simplifies our description to assume that $b$ is input by a $(d\times d)$-matrix $B=sB^{\sigma t}$ with $s=\pm 1$ and $\sigma^2=1$.
If the $s$ and $\sigma$ are not specified with $B$, then suitable values can be detected during the execution of the algorithms for \thmref{thm:main}.  We either prove that $B=0$ or we find $u,v\in V$ such that $b(u,v)=uBv^{\sigma t}\neq 0$.  The first such pair $u,v\in V$ determines $\sigma$ by $\sigma:\alpha\mapsto b(u,\alpha v)b(u,v)^{-1}$, and $s=b(v,u)^{-1}\cdot b(u,v)^{\sigma}$.   We write the algorithm as though $s$ and $\sigma$ are known.

%----------------------------
\section{Smiley's method}\label{sec:cubic}

Let us start with the algorithm {\tt Decompose} which is not asymptotically optimal, but which (I believe) is the simplest to implement and captures all Hermitian forms at once.  This is the prototype for the optimal sequential and parallel algorithms given later. 

If $A$ is an elementary matrix then $ABA^{\sigma t}$ modifies $B$ in one of three ways.  First, if $A$ is a diagonal matrix with $1$'s on the diagonal except for $\lambda$ in entry $i$, then $ABA^{\sigma t}$ scales row $i$ by $\lambda$, and column $i$ by $\lambda^{\sigma}$.  For instance:
\begin{equation*}
	\begin{bmatrix} 1 &  & \\  & \lambda & \\ &  & 1 \end{bmatrix}
	\begin{bmatrix} \alpha & \beta & \delta \\ s\beta^{\sigma} & \gamma & \epsilon\\ 
		s\delta^{\sigma} & s\epsilon^{\sigma} & \phi
	\end{bmatrix}
	\begin{bmatrix} 1 &  & \\  & \lambda & \\ &  & 1 \end{bmatrix}^{\sigma t}
	=
	\begin{bmatrix} \alpha & \beta\lambda^{\sigma} & \delta \\ 
		s(\beta\lambda^{\sigma})^{\sigma} & \lambda \gamma\lambda^{\sigma} & \lambda\epsilon\\ 
		s\delta^{\sigma} & s(\lambda\epsilon)^{\sigma} & \phi
	\end{bmatrix}
\end{equation*}
We describe that as \emph{scaling row-column $i$ by $\lambda$}.
Second, if $A$ is a transposition of $i$ and $j$ then $ABA^{\sigma t}$ has the entries from $B$ with rows $i$ and $j$ swapped as well as columns $i$ and $j$ swapped.  We call this \emph{swapping row-column $i$ with row-column $j$}.  That does not involve operations in $\Delta$.  Thirdly, if $A=I+\lambda E_{ij}$ then $ABA^{\sigma t}$ has the effect of adding $\lambda$ times row $i$ to row $j$ and $\lambda^{\sigma}$ times column $i$ to column $j$, as illustrated below.
\begin{equation*}
	\begin{bmatrix} 1 &  & \\  & 1 & \\ & -\gamma^{\sigma} & 1 \end{bmatrix}
	\begin{bmatrix} 0 & 1 & \gamma \\ s & \beta & \delta\\ 
		s\gamma^{\sigma} & s\delta^{\sigma} & \epsilon
	\end{bmatrix}
	\begin{bmatrix} 1 &  & \\  & 1 & \\ & -\gamma^{\sigma} & 1 \end{bmatrix}^{\sigma t}
	=
	\begin{bmatrix} 0 & 1 & 0 \\ 
		s & \beta & \delta -\beta\gamma\\ 
		0 & s(\delta-\beta\gamma)^{\sigma} & *
	\end{bmatrix}
\end{equation*}
That implicitly involved the fact that entries $\beta$ on the diagonal satisfy $\beta=s\beta^{\sigma}$.  To \emph{clear a row-column} means to use a selected non-zero entry $j$ in a row-column $i$ of $B$, and use successive multiplications by $I+\lambda_k E_{ki}$, for $k\in \{1,\dots,d\}-\{i\}$ to set all other entries in the row-column $i$ to zero. This is possible whenever $i=j$ or $B_{ii}=0$.  Using the symmetry of the matrices $B=sB^{\sigma t}$, clearing a row-column uses $d^2+O(d)$ additions, $d^2+O(d)$ multiplications, $d$ applications of $\sigma$, and one inversion.

We use upper case Roman letters for block sub-matrices and lower case Greek letters for coefficients in $\Delta$.  We also assume that the associated matrix $A$ which transforms $B$ into the return $ABA^{\sigma t}$ as in \eqref{eq:standard} is evident form the operations described, and so we do not explicitly include $A$ in the description of the algorithm.

%------------ STANDARD BLOCK
\begin{figure}[!htbp]
%\rule{\textwidth}{1pt}
\flushleft\texttt{Standardize$\left(~B=\begin{bmatrix} 0 & 1 \\ s & \alpha\end{bmatrix}~\right)$:}\\
\begin{minipage}{0.9\textwidth}
If $\alpha\neq 0$, set $A=\begin{bmatrix} 1 & -\alpha^{-1} \\ 0 & 1 \end{bmatrix}$ and
return $ABA^{\sigma t}=[-s\alpha^{-1}]\oplus[\alpha]$. If $\alpha=0$ and $s=1\neq -1$, set $A=\begin{bmatrix} 1 & 1 \\ 1 & -1 \end{bmatrix}$ and return $ABA^{\sigma t}=[2]\oplus [-2]$.
Else return $B$.
\end{minipage}
%\rule{\textwidth}{1pt}
\end{figure}
%-----------

%------------ STANDARD BASIS
\begin{figure}[!htbp]
%\rule{\textwidth}{1pt}
\flushleft\texttt{Decompose$(~B\in M_d(\Delta) : B=sB^{\sigma t}~)$:}
\begin{enumerate}[(I)]
\item \textbf{[Base case]} If $d\leq 1$ return $B$.

\item\label{A.step:ansi} \textbf{[Anisotropic case]} If $B_{11}=\beta\neq 0$, 
then use that entry to clear the remaining non-zero entries of row-column 1.  
Now $B=\begin{bmatrix} \beta & 0 \\ 0 & B'\end{bmatrix}$.\\
\noindent Return $[\beta]\oplus {\tt Decompose}(B')$.

\item\label{A.step:iso} \textbf{[Isotropic case]} Else, if $B_{12}=\gamma\neq 0$ (after a possible swap of
a row-column), i.e. $B=\begin{bmatrix} 0 & \gamma & * \\ s \gamma^{\sigma} & \alpha & * \\ * & * & * \end{bmatrix}$, then scale row-column $2$ by $\gamma^{-1}$ and excluding $B_{22}$, use $B_{12}$ to clear row-column $1$  and $B_{21}$ to clear row-column $2$.  Now
$B=\begin{bmatrix} B' & 0 \\ 0 & B'' \end{bmatrix}$ where $B'=\begin{bmatrix} 0 & 1\\ s & \alpha\end{bmatrix}$.\\ 
\noindent Return ${\tt Standardize}(B')\oplus {\tt Decompose}(B'')$.

\item\label{A.step:rad} \textbf{[Radical case]} Else, 
$B=\begin{bmatrix} 0 & 0 \\ 0 & B' \end{bmatrix}$ so return $[0]\oplus {\tt Decompose}(B')$.
\end{enumerate}
%\rule{\textwidth}{1pt}
\end{figure}
%-----------

\begin{proof}[Proof of \thmref{thm:main}(i)]
The algorithm \texttt{Decompose} returns a block diagonal matrix whose blocks are as in \eqref{eq:standard}.    That algorithm only modifies the entries of $B$ so that the space complexity is $O(d^2)$ elements in $\Delta$.  

Now we consider the time complexity.  There are at most $d$ equality tests to decide on the correct case to enter.  The anisotropic case clears one row-column and recurses on a matrix of dimension $d-1$.  The isotropic case clears two row-columns, performs some multiplications of $(2\times 2)$-matrices, and recurses on a matrix of dimension $d-2$.  Finally, the radical case simply recurses on a matrix of dimension $d-1$.  Hence, if $T(d)$ is the number of additions performed by the algorithm, then $T(d)\in 2d^2+T(d-2)+O(d)$.  If $r$ is the dimension of the radical and $e=d-r$, then $T(d)\in e^3/3+O(d^2)$.  The algorithm uses the same number of multiplications, $\binom{d}{2}$ equality tests, $e$ inversions, and $\binom{d}{2}-\binom{r}{2}$ applications of $\sigma$.  
\end{proof}

%----------------------------------------------------
\section{Optimal and parallel methods}\label{sec:optimal}

Multiplication of $(d\times d)$-matrices by the traditional algorithm is not the most efficient method for large dimensions.  The various new methods use $O(d^{\omega})$ operations in $\Delta$ for some $2\leq \omega\leq 3$ \cite[p. 315]{vzG}.  Here we prove the same complexity for finding a decomposition as in \eqref{eq:standard}.  B\"{u}rgisser et. al. give an example of a symmetric $(d\times d)$-matrix over a field where the complexity of finding an orthogonal basis is $O(d^{\omega})$ (provided that $\omega>2$) \cite[Theorem 16.20]{Burgisser:comp} and so the complexity in \thmref{thm:main}(ii) is best possible
in general.  

\begin{figure}[!htbp]
%\rule{\textwidth}{1pt}
\flushleft\texttt{DecomposeByBlocks$(~B\in M_d(\Delta) : B=sB^{\sigma t}~)$:}
\begin{enumerate}[(I)]
\item\label{step:radical} \textbf{[Detect Radical]}
Compute an invertible $A$ such that $AB=\begin{bmatrix} B'\\ 0 \end{bmatrix}$ where $B'$ has full row rank.  Now $ABA^{\sigma t}=\begin{bmatrix} B'' & 0 \\ 0 & 0 \end{bmatrix}$ with $B''$ nonsingular.  Apply step \eqref{step:anisotropic} to $B''$.

\item\label{step:anisotropic} \textbf{[Block Anisotropic case]}
Here $B$ is a nonsingular $(d\times d)$-matrix.  If $d\leq 1$, halt; else take $B=\begin{bmatrix} B' & * \\ * & * \end{bmatrix}$ with $B'\in M_{\lceil d/2\rceil}(\Delta)$.   Find $A$ such that $AB'A^{\sigma t} = \begin{bmatrix} B'' & 0 \\ 0 & 0 \end{bmatrix}$ and $B''$ has full rank (as in (I)).  
Compute
\begin{equation}\label{eq:ani-1}
	 \begin{bmatrix} A & 0 \\ 0 & I \end{bmatrix}
	\begin{bmatrix} B' & * \\ * & * \end{bmatrix}
	\begin{bmatrix} A & 0 \\ 0 & I \end{bmatrix}^{\sigma t}
	=\begin{bmatrix} AB' A^{\sigma t} & * \\ * & *\end{bmatrix}
	=\begin{bmatrix} 
			B'' & 0 & C \\
			0 & 0 & W \\
			sC^{\sigma t} & sW^{\sigma t} & *
	\end{bmatrix}.
\end{equation}
If $B''$ has dimension $0$ then (as $B$ is nonsingular) $W$ is nonsingular; proceed to step \eqref{step:Hyperbolic}.  Otherwise, $B''=s(B'')^{\sigma t}$ is nonsingular.  Set $Y=-sC^{\sigma t}(B'')^{-1}$ and compute
\begin{equation}\label{eq:ani-2}
\begin{bmatrix} 
	I & 0 & 0 \\
	0 & I & 0 \\
	Y & 0 & I
\end{bmatrix}	
\begin{bmatrix} 
	B'' & 0 & C \\
	0 & 0 & * \\
	sC^{\sigma t} & * & *
\end{bmatrix}
\begin{bmatrix} 
	I & 0 & Y^{\sigma t} \\
	0 & I & 0 \\
	0 & 0 & I
\end{bmatrix}
=
\begin{bmatrix}
	B'' & 0 & 0 \\
	0 & 0 & X \\
	0 & sX^{\sigma t} & Z
\end{bmatrix}.
\end{equation}
Note $X$ has full row rank since $B$ is nonsingular.
Apply step \eqref{step:anisotropic} to $B''$, and apply step \eqref{step:Hyperbolic} to $\begin{bmatrix} 0 & X \\ sX^{\sigma t} & Z\end{bmatrix}$; then halt.

\item\label{step:Hyperbolic} \textbf{[Block Isotropic case]}
Now $B=\begin{bmatrix} 0 & X \\ sX^{\sigma t} & *\end{bmatrix}$ and $X$ has full row rank.  Compute an invertible matrix $A$ such that $XA=\begin{bmatrix} C &  0\end{bmatrix}$ where $C$ has full column rank; thus, $C$ is invertible.  Compute
\begin{equation}\label{eq:hype-1}
\begin{bmatrix} C^{-1} & 0 \\ 0 & A^{\sigma t}\end{bmatrix}
\begin{bmatrix} 0 & X \\ sX^{\sigma t} & *\end{bmatrix}
\begin{bmatrix} C^{-\sigma t} & 0 \\ 0 & A\end{bmatrix}
=
\begin{bmatrix}
	0 & I & 0 \\
	sI & Z & Y\\
	0 & sY^{\sigma t} & B'
\end{bmatrix}.
\end{equation}
Observe that:
\begin{equation}\label{eq:hype-2}
\begin{bmatrix}
	I & 0 & 0 \\
	0 & I & 0 \\
	-sY^{\sigma t} & 0 & I
\end{bmatrix}
\begin{bmatrix}
	0 & I & 0 \\
	sI & Z & Y\\
	0 & sY^{\sigma t} & B'
\end{bmatrix}
\begin{bmatrix}
	I & 0 & -sY \\
	0 & I & 0 \\
	0 & 0 & I
\end{bmatrix}
=
\begin{bmatrix}
	0 & I & 0 \\
	sI & Z & 0\\
	0 & 0 & B'
\end{bmatrix}.
\end{equation}
Let $B''=\begin{bmatrix} 0 & I\\ sI & Z\end{bmatrix}$ and decompose 
$Z=s Z^{\sigma t}=U+D+sU^{\sigma t}$ where $U$ is upper triangular with $0$ entries on the diagonal.  So $D$ is diagonal and $D=sD^{\sigma}$.  Reset $B''$ to be
\begin{equation}\label{eq:hype-3}
\begin{bmatrix}
	I & 0 \\
	-U & I
\end{bmatrix}
\begin{bmatrix}
	0 & I \\
	sI & Z 
\end{bmatrix}
\begin{bmatrix}
	I & -U^{\sigma t} \\
	0 & I 
\end{bmatrix}
=
\begin{bmatrix}
 0 & I \\
sI & D
\end{bmatrix}.
\end{equation}
Sort the row-columns so that the matrix is in the form 
$\begin{bmatrix} 0 & 1 \\ s & \alpha_1\end{bmatrix}\oplus \cdots\oplus\begin{bmatrix} 0 & 1\\
s & \alpha_f\end{bmatrix}$ with $\alpha_i\in \Delta$.  Apply ${\tt Standardize}$ to each of those blocks.  Finally, $B'$ is nonsingular so apply step \eqref{step:anisotropic} to $B'$, then halt.
\end{enumerate}
%\rule{\textwidth}{1pt}
\end{figure}

\begin{proof}[Proof of \thmref{thm:main}(ii)]
The algorithm {\tt DecomposeByBlocks} suffices as described so it remains to analyze the time complexity of the algorithm.

 We start by detecting the radical of $B$.  This amounts to solving for a basis of the null space of $B$.  That has complexity of $O(d^{\omega})$ \cite[Theorem 2]{Sol}.  To create $B''$ requires 2 matrix multiplications and $d^2$ applications of $\sigma$.  Thus the radical case uses $O(d^{\omega})$ operations in $\Delta$.  The algorithm never re-enters this case.
 
In the block anisotropic case we solve for a null space on a $\lceil d/2\rceil$-square matrix $B'$, and multiply two $(d\times d)$-matrices, in \eqref{eq:ani-1}.
Let $f$ be the rank of the null space of $B'$.
At this point we have two cases.  If $f=\lceil d/2\rceil$ we exit the block anisotropic case and enter the block semi-hyperbolic case; otherwise, we to create $Y$ (we invert and multiply a $\big((\lceil d/2\rceil-f)\times (\lceil d/2\rceil-f)\big)$-matrix), multiply two $(d\times d)$-matrices in \eqref{eq:ani-2}.  We then make one recurse call to the block anisotropic case for $B''$, and one call to the block semi-hyperbolic case for $\begin{bmatrix} 0 & X\\ sX^{\sigma t} & Z\end{bmatrix}$ where $X$ has $f$ rows
and $Z$ is $(d-\lceil d/2\rceil)\times (d-\lceil d/2\rceil)$.  Ignoring the recursions, the anisotropic case uses $O(d^{\omega})$ operations in $\Delta$.  

The block semi-hyperbolic case takes in a $(d\times d)$-matrix partitioned by into $(f,d-f)$-blocks. We compute a null column space of an $(f\times (d-f))$-matrix, multiply 2 $(d\times d)$-matrices \eqref{eq:hype-1} (\eqref{eq:hype-2} requires no computation), and we also 
multiply two $(2f\times 2f)$-matrices in \eqref{eq:hype-3}.  Finally there are at most $d/2$ applications of ${\tt Standardize}$ and a recursive call on a $((d-2f)\times (d-2f))$-matrix $B'$.  All this amounts to $O(d^{\omega})$ operations in $\Delta$ before the recursion.

Now we estimate the total cost.
Let $T_{a}(d)$ be the cost of the block anisotropic case for an input of dimension $d$, 
and $T_{h}(d,f)$ the cost of the block semi-hyperbolic case for an input of dimension $d$
where $X$ has $f$-rows.  For some constants $C_a,C_h>0$,
$T_{h}(d,f) \leq  T_a(d-2f)+C_h d^{\omega}$, and
\begin{align*} 
	T_{a}(d) & \leq \max\{~T_{h}(d,d/2),~T_{a}(d/2-f)+T_{h}(d/2+f,f))~\}+C_a d^{\omega}\\
	 	& \leq 2T_a(d/2-f)+C_h (d/2+f)^{\omega}+(C_a+C_h) d^{\omega}\\
		& \leq 2T_a(d/2) + (C_a+2C_h) d^{\omega}.
\end{align*}
Thus, $T(d)\in O(d^{\omega})$.
\end{proof}

\begin{proof}[Proof of \thmref{thm:main}(iii)] {\tt DecomposeByBlocks} uses $O(\log d)$ recursive calls and each step can use the parallel $NC^2$ (i.e. $O(\log^2 d)$) linear algebra algorithms of \cite[Sections 3.8, 4.5]{KR:parallel} and \cite[Section 2.3]{Sol} to find null-spaces and multiply matrices in an arithmetic model.  (Those methods make it possible to trade on time efficiency to reduce the number of required processors, which is of importance in practice.)
\end{proof}

%==============================================
\section{Post-processing adjustments}\label{sec:adjust}

In our algorithms we opted for a decomposition of $B$ which is as close to diagonal as possible so that the associated basis is nearly orthogonal.  It is also common to want a decomposition with as many blocks of the form $J=\begin{bmatrix} 0 & 1\\ s & 0 \end{bmatrix}$ as possible.  The algorithm can be tuned in that direction by modifying {\tt Standardize} and by converting various $(1\times 1)$-blocks into $J$'s at the end of the algorithm.  The details are analogous to those used in {\tt Standardize}.

In some cases a canonical return is possible with a few adjustments.  For example,   
the block $[\alpha]$ can be adjusted to $[\gamma \alpha \gamma^{\sigma}]$ for $0\neq \gamma\in \Delta$.  Hence, if $\alpha=\gamma^{-1}\gamma^{-\sigma}$ for some $\gamma\in \Delta-\{0\}$, then we may replace $\alpha$ with $1$. Computationally finding $\gamma$ to perform this adjustment can be involved.  Already when $\sigma=1$ and $\Delta$ a field this amounts to finding a square-root of $\alpha$.  If the number of classes in $\Delta$ of the form $\gamma\alpha\gamma^{\sigma}$ is linearly ordered, it is possible to sort the $(1\times 1)$-blocks accordingly.

Another situation for modification tries to convert multiple $(1\times 1)$-blocks.
For instance, if $\Delta$ is a field and $0\neq \alpha = \gamma\gamma^{\sigma}+\delta\delta^{\sigma}$ for some $\gamma,\delta\in \Delta$ (for example, if $\sigma=1$ and $\alpha$ is a sum of squares), then 
\begin{equation*}
	\begin{bmatrix} \gamma & \delta \\ \delta^{\sigma} & -\gamma^{\sigma} \end{bmatrix}
	\begin{bmatrix} 1 & 0 \\ 0 & 1  \end{bmatrix}
	\begin{bmatrix} \gamma & \delta \\ \delta^{\sigma} & -\gamma^{\sigma} \end{bmatrix}^{\sigma t}
=\begin{bmatrix} \alpha & 0 \\ 0 & \alpha \end{bmatrix}.
\end{equation*} 
Similarly, if the characteristic is $2$ and $\alpha\neq 0$ then
\begin{equation*}
\begin{bmatrix} 0 & \alpha & 1\\ 1 & \alpha & 1 \\ 1 & 0 & 1\end{bmatrix}
\begin{bmatrix} 0 & 1 &   \\ 1 & 0 &  \\  &  & \alpha \end{bmatrix}
\begin{bmatrix} 0 & \alpha & 1\\ 1  & \alpha & 1 \\ 1 & 0 & 1\end{bmatrix}^{\sigma t}
=
\begin{bmatrix} \alpha &  & \\  & \alpha &  \\  &  & \alpha \end{bmatrix}.
\end{equation*}

\section*{Acknowledgements}
Thanks to Peter Brooksbank for suggesting this note and offering comments.

\bibliographystyle{amsplain}

\end{document}